\documentclass[12pt]{article}
\usepackage{amsfonts}

\begin{document}

\renewcommand{\clubsuit}{\quad \vrule height7.5pt width4.17pt depth0pt}

\newcommand {\IN}{\mathbf{N}}
\newcommand {\IR}{\mathbf{R}}
\newcommand {\IQ}{\mathbf{Q}}

\newcommand {\D}{\mbox{${\mathcal D}$}}
\newcommand {\E}{\mbox{${\mathcal E}$}}
\newcommand {\F}{\mbox{${\mathcal F}$}}
\newcommand {\G}{\mbox{${\mathcal G}$}}
\newcommand {\IH}{\mbox{${\mathcal H}$}}

\newcommand {\IP}{\mbox{${\mathcal P}$}}

\newcommand {\IS}{\mbox{${\mathcal S}$}}

\newcommand{\half}{\frac{1}{2}}
\newcommand{\proof}{\noindent {\sc Proof.\ \ }}  
\newcommand{\ti}[1]{\tilde{#1}}

\newtheorem{stat}{Statement}
  \newtheorem{examp}[stat]{Example}
  \newtheorem{assump}{Assumption}[section]
  \newtheorem{decth}[stat]{Theorem}
  \newtheorem{prop}[stat]{Proposition}
  \newtheorem{cor}[stat]{Corollary}
  \newtheorem{thm}[stat]{Theorem}
  \newtheorem{lemma}[stat]{Lemma}
  \newtheorem{remark}[stat]{Remark}
  \newtheorem{def1}{Definition}[section]

\title{Some non-linear s.p.d.e.'s that are second order in time}
\author{Robert C. Dalang$^1$ and Carl Mueller$^2$
\\
\\D\'epartement de Math\'ematiques
\\Ecole Polytechnique F\'ed\'erale
\\1015 Lausanne, Switzerland
\\E-mail:  robert.dalang@epfl.ch
\\
\\Department of Mathematics
\\University of Rochester
\\Rochester, NY  14627  USA
\\E-mail:  cmlr@troi.cc.rochester.edu}
\date{}
\maketitle

\footnotetext[1]{Supported in part by the
Swiss National Foundation for Scientific Research.}
\footnotetext[2]{Supported by an NSA grant. 


AMS 2000 {\em subject classifications}
Primary, 60H15; Secondary, 35R60, 35L05.}

\begin{abstract} 
We extend Walsh's theory of martingale measures in order to deal with 
hyperbolic stochastic partial differential equations that are second order in
time, such as the wave 
equation and the beam equation, and driven by spatially homogeneous Gaussian noise.  For such equations, the fundamental 
solution can be a distribution in the sense of Schwartz, which appears as an
integrand in the reformulation of the s.p.d.e.~as a stochastic integral
equation.  Our approach 
provides an alternative to the Hilbert space integrals of Hilbert-Schmidt
operators.  We give several examples, including the beam equation and 
the wave equation, with nonlinear multiplicative noise terms.  
\end{abstract}

\section{Introduction}
\setcounter{equation}{0}
The study of stochastic partial differential equations (s.p.d.e.'s) began 
in earnest following the papers of Pardoux \cite{par72}, \cite{par75}, 
\cite{par77}, and Krylov and Rozovskii \cite{kr81}, \cite{kr82}. 
Much of the literature has been concerned with the heat equation, most often 
driven by space-time white noise, and with related parabolic equations.  
Such equations are first order in time, and generally second order in the
space variables.  There has been much less work on s.p.d.e.'s that are second
order in time, such as the wave equation and related hyperbolic equations.
Some early references are Walsh \cite{wal86}, and Carmona and Nualart
\cite{cn88a}, \cite{cn88b}.  More recent papers are Mueller \cite{mue97}, 
Dalang and Frangos \cite{df98}, and Millet and Sanz-Sol\'e \cite{mill99}.

For linear equations, the noise process can be considered as a random 
Schwartz distribution, and therefore the theory of deterministic p.d.e.'s 
can in principle be used.  However, this yields solutions in the space of 
Schwartz distributions, rather than in the space of function-valued stochastic 
processes.  For linear s.p.d.e.'s such as the heat and wave equation driven 
by space-time white noise, this situation is satisfactory, since, in fact, 
there is no function-valued solution when the spatial dimension is greater
than $1$.  However, since non-linear functions of Schwartz distributions are difficult 
to define (see however Oberguggenberger and Russo \cite{or98}), it is difficult to 
make sense of non-linear s.p.d.e.'s driven by space-time white noise in 
dimensions greater than $1$.

A reasonable alternative to space-time white noise is Gaussian noise with some 
spatial correlation, that remains white in time.  This approach has been 
taken by several authors, and the general framework is given in Da Prato 
and Zabczyk \cite{dz92}.  However, there is again a difference 
between parabolic and hyperbolic equations: while the Green's function is 
smooth for the former, for the latter it is less and less regular as the 
dimension increases.  For instance, for the wave equation, the Green's 
function is a bounded function in dimension 1, is an unbounded function in 
dimension 2, is a measure in dimension 3, and a Schwartz distribution in 
dimensions greater than 3.

There are at least two approaches to this issue.  One is to extend the 
theory of stochastic integrals with respect to martingale measures, as 
developed by Walsh \cite{wal86}, to a more general class of integrands that 
includes distributions.  This approach was taken by Dalang \cite{dal99}.  In the 
case of the wave equation, this yields a solution to the non-linear 
equation in dimensions 1, 2 and 3.  The solution is a {\em random field,} 
that is, it is defined for every $(t,x) \in \IR_+ \times \IR^d$.  Another 
approach is to consider solutions with values in a function space, 
typically an $L^2$-space: for each fixed $t \in \IR_+$, the solution is an 
$L^2$-function, defined for almost all $x \in \IR^d$.  This approach has 
been taken by Peszat and Zabczyk in \cite{pz00} and Peszat \cite{pes02}. 
In the case of the non-linear wave equation, this approach yields
function-valued solutions in all dimensions. It should be noted that the notions
of {\em random field} solution and {\em function-valued} solution are {\em not}
equivalent: see L\'ev\`eque \cite{Lev01}.

In this paper, we develop a general approach to non-linear s.p.d.e's, with 
a focus on equations that are second order in time, such as the wave 
equation and the beam equation.  This approach goes in the direction of 
unifying the two described above, since we begin in Section \ref{sec2} with 
an extension of Walsh's martingale measure stochastic integral \cite{wal86}, 
in such a way as to integrate processes that take values in an $L^2$-space, 
with an integral that takes values in the same space.  This extension defines 
stochastic integrals of the form
\[
   \int_0^t\int_{\IR^d} G(s,\cdot-y) Z(s,y)\, M(ds,dy),
\]
where $G$ (typically a Green's function) takes values in the space of 
Schwartz distributions, $Z$ is an adapted process with values in 
$L^2(\IR^d)$, and $M$ is a Gaussian martingale measure with spatially 
homogeneous covariance.

With this extended stochastic integral, we can study non-linear forms of a 
wide class of s.p.d.e.'s, that includes the wave and beam equations in all 
dimensions, namely equations for which the p.d.e. operator is
\[
   \frac{\partial^2 u}{\partial t^2} + (-1)^k \Delta^{k} u,
\]
where $k \geq 1$ (see Section \ref{sec3}).  Indeed, in Section \ref{sec4} 
we study the corresponding non-linear s.p.d.e.'s.  We only impose the minimal 
assumptions on the spatial covariance of the noise, that are needed even for 
the linear form of the s.p.d.e. to have a function-valued solution.  The 
non-linear coefficients must be Lipschitz and vanish at the origin.  This 
last property guarantees that with an initial condition that is in 
$L^2(\IR^d)$, the solution remains in $L^2(\IR^d)$ for all time.

In Section \ref{sec5}, we specialize to the wave equation in a weighted 
$L^2$-space, and remove the condition that the non-linearity vanishes at 
the origin.  Here, the compact support property of the Green's function of 
the wave equation is used explicitly.  We note that Peszat \cite{pes02} also 
uses weighted $L^2$ spaces, but with a weight that decays exponentially at 
infinity, whereas here, the weight has polynomial decay.

\section{Extensions of the stochastic integral}\label{sec2}
\setcounter{equation}{0}
In this section, we define the class of Gaussian noises that drive the 
s.p.d.e.'s that we consider, and give our extension of the martingale 
measure stochastic integral.

Let ${\cal{D}}(\IR^{d+1})$ be the topological vector space of functions 
$\varphi \in C^\infty_0(\IR^{d+1})$, the space of infinitely 
differentiable functions with compact support, with the standard notion of 
convergence on this space (see Adams \cite{ada75}, page 19). 
Let $\Gamma$ be a non-negative and {\em non-negative definite} 
(therefore symmetric) tempered measure on $\IR^d$.  That is,
\[
   \int_{\IR^d} \Gamma (dx)\, (\varphi \ast \tilde{\varphi})(x) \geq 0, 
    \qquad \mbox{for all } \varphi \in {\cal{D}}(\IR^d),
\]
where $\tilde{\varphi}(x) = \varphi(-x)$, ``$\ast$" denotes convolution, 
and there exists $r > 0$ such that
\begin{equation}\label{2.1}
   \int_{\IR^d} \Gamma(dx)\, \frac{1}{(1+\vert x \vert^2)^r} < \infty.
\end{equation}
We note that if $\Gamma(dx) = f(x) dx$, then
\[
   \int_{\IR^d} \Gamma(dx)(\varphi \ast \tilde{\varphi})(x) 
   = \int_{\IR^d} dx \int_{\IR^d} dy \,       \varphi (x) f(x-y) \varphi(y)
\]
(this was the framework considered in Dalang \cite{dal99}).  Let $\IS(\IR^d)$ 
denote the Schwartz space of rapidly decreasing $C^\infty$ test functions, and 
for $\varphi \in \IS(\IR^d)$, let $\F\varphi$ denote the Fourier transform 
of $\varphi$:
\[
   \F\varphi(\eta) = \int_{\IR^d} \exp(-i\, \eta \cdot x) \varphi(x)\, dx.
\]
According to the Bochner-Schwartz theorem (see Schwartz \cite{sch66}, 
Chapter VII, Th\'eor\`eme XVII), there is a non-negative tempered measure 
$\mu$ on $\IR^d$ such that $\Gamma = \F\mu,$ that is
\begin{equation}\label{2.2}
   \int_{\IR^d} \Gamma (dx)\, \varphi(x) = \int_{\IR^d} \mu (d\eta)\, \F\varphi(\eta), 
    \qquad \mbox{for all } \varphi \in \IS(\IR^d).
\end{equation}

\noindent{\bf Examples.} (a) Let $\delta_0$ denote the Dirac functional. 
Then $\Gamma(dx) = \delta_0 (x)\, dx$ satisfies the conditions above. 

(b) Let $0 < \alpha < d$ and set $f_\alpha(x) = \vert x \vert^{-\alpha}$, 
$x \in \IR^d$.  Then $f_\alpha = c_\alpha \F f_{d-\alpha}$ 
(see Stein \cite{ste70}, Chapter V \S 1, Lemma 2(a)), 
so $\Gamma(dx) = f_\alpha(x)\, dx$ also satisfies the conditions above.
\vskip 16pt

Let $F = (F(\varphi),\ \varphi \in {\cal{D}}(\IR^{d+1}))$ be an 
$L^2(\Omega, \G, P)$-valued mean zero Gaussian process with covariance 
functional
\[
   (\varphi, \psi) \mapsto E(F(\varphi) F(\psi)) 
   = \int_{\IR^d} \Gamma(dx)\, (\varphi \ast        \tilde{\psi})(x).
\]
As in Dalang and Frangos \cite{df98} and Dalang \cite{dal99}, 
$\varphi \mapsto F(\varphi)$ extends to a worthy martingale measure 
$(t, A) \mapsto M_t(A)$ (in the sense of Walsh \cite{wal86}, pages 289-290) 
with covariance measure
\[
   Q([0, t] \times A \times B) = \langle M(A), M(B)\rangle_t 
   = t \int_{\IR^d} \Gamma (dx)      \int_{\IR^d} dy\, 1_A(y)\, 1_B(x+y)
\]
and dominating measure $K \equiv Q,$ such that
\[
   F(\varphi) = \int_{\IR_+} \int_{\IR^d} \varphi(t, x) M(dt, dx), 
   \qquad \mbox{for all } \varphi       \in {\cal{D}}(\IR^{d+1}).
\]
The underlying filtration is $(\F_t = \F_t^0 \vee {\cal{N}},\ t \geq 0)$, where 
\[
   \F^0_t= \sigma (M_s(A), \ s \leq t, \ A \in {\cal{B}}_b(\IR^d)),
\]
${\cal{N}}$ is the $\sigma$-field generated by $P$-null sets and 
${\cal{B}}_b(\IR^d)$ denotes the bounded Borel subsets of $\IR^d$.

Recall \cite{wal86} that a function 
$(s, x, \omega) \mapsto g(s, x; \omega)$ is termed {\em elementary} if it 
is of the form
\[
   g(s, x; \omega) = 1_{]a, b]}(s) 1_A(x) X(\omega),
\]
where $0 \leq a < b$, $A \in {\cal{B}}_b(\IR^d)$ and $X$ is a bounded and 
$\F_a$-measurable random variable.  The $\sigma$-field on 
$\IR_+ \times \IR^d \times \Omega$ generated by elementary functions is 
termed the {\em predictable} $\sigma$-field.

Fix $T > 0$.  Let ${\cal{P}}_+$ denote the set of predictable functions 
$(s, x; \omega) \mapsto g(s, x; \omega)$ such that 
$\Vert g \Vert_+ < \infty$, where
\[
   \Vert g \Vert^2_+ = E\left(\int_0^T ds \int_{\IR^d} \Gamma(dx) 
   \int_{I\!\!R^d} dy\, \vert g (s, y) g(s, x + y)\vert\right).
\]
Recall \cite{wal86} that ${\cal{P}}_+$ is the completion of the set of 
elementary functions for the norm $\Vert \cdot \Vert_+$. 

For $g \in {\cal{P}}_+$, Walsh's stochastic integral
\[
   M_t^g(A) = \int_0^t ds \int_A g(s, x) M(ds,dx)
\]
is well-defined and is a worthy martingale measure with covariation measure
\[
   Q_g([0, t] \times A \times B) 
   = \int_0^t ds \int_{\IR^d} \Gamma (dx) \int_{\IR^d} dy\, 1_A(y) 1_B(x+y)\, g(s, y) g(s, x+y)
\]
and dominating measure
\[
   K_g([0, t] \times A\times B) = \int_0^t ds \int_{\IR^d} \Gamma(dx) 
   \int_{\IR^d} dy\, 1_A(y) 1_B(x + y)\, \vert g(x, y) g(s, x+y)\vert.
\]

For a deterministic real-valued function $(s, x) \mapsto g(s, x)$ and a 
real-valued stochastic process $(Z(t, x),\ (t, x) \in \IR_+ \times \IR^d)$, 
consider the following hypotheses ($T>0$ is fixed).

\begin{itemize}
\item[{\bf (G1)}] For $0 \leq s \leq T$, $g(s,\cdot) \in C^\infty(\IR^d)$, 
$g(s,\cdot)$ is bounded uniformly in $s$, and $\F g(s,\cdot)$ is a function.

\item[{\bf (G2)}] For $0 \leq s \leq T$, $Z(s, \cdot) \in C_0^\infty(\IR^d)$ 
a.s., $Z(s, \cdot)$ is $\F_s$-measurable, and in addition, there is a 
compact set $K \subset \IR^d$ such that supp~$Z(s,\cdot) \subset K$, for 
$0 \leq s \leq T$.  Further, $s \mapsto Z(s,\cdot)$ is mean-square 
continuous from $[0,T]$ into $L^2(\IR^d)$, that is, for $s \in [0,T]$,
\[
   \lim_{t \to s} E\left(\Vert Z(t,\cdot) - Z(s,\cdot) \Vert^2_{L^2(\IR^d)}\right) = 0.
\]

\item[{\bf (G3)}] $I_{g, Z} < \infty$, where
\[
  I_{g, Z} = \int_0^T ds \int_{\IR^d} d \xi\, E(\vert \F Z(s, \cdot)(\xi)\vert^2) \int_{\IR^d} 
      \mu(d \eta)\, \vert \F g(s, \cdot)(\xi-\eta)\vert^2.
\]
\end{itemize}

\begin{lemma} Under hypotheses (G1), (G2) and (G3), for all $x \in \IR^d$, 
the function defined by $(s, y; \omega) \mapsto g(s, x-y) Z(s, y; \omega)$ 
belongs to ${\cal{P}}_+$, and so
\[
  v_{g, Z} (x) = \int^T_0 \int_{\IR^d} g(s, x-y) Z(s, y) M(ds, dy)
\]
is well-defined as a (Walsh-) stochastic integral.  Further, a.s., 
$x \mapsto v_{g,Z}(x)$ belongs to $L^2(\IR^d)$, and
\begin{equation}\label{2.3}
   E\left(\Vert v_{g,Z} \Vert_{L^2(\IR^d)}^2\right) = I_{g,Z}.
\end{equation}
\label{rlem1}
\end{lemma}

\proof Observe that $\Vert g(\cdot, x-\cdot) Z(\cdot, \cdot)\Vert_+^2$ is equal to
\[
   E\left(\int_0^T ds \int_{\IR^d} \Gamma (dz) \int_{\IR^d} dy\, 
   \vert g(s, x-y) Z(s, y) g(s, x - y-z) Z(s, y + z)\vert\right).
\]
Because $g(s,\cdot)$ is bounded uniformly in $s$ by (G1), this expression is bounded by a constant times
\begin{eqnarray*}
  && E\left(\int_0^T ds \int_{\IR^d} \Gamma (dz) \int_{\IR^d} dy\, \vert Z(s, y) Z(s, y + z)\vert\right) \\
  &&\qquad   = E\left(\int_0^T ds \int_{\IR^d} \Gamma (dz)\, 
  (\vert Z(s,\cdot)\vert \ast \vert    \tilde{Z}(s,\cdot)\vert)(-z)  \right).
\end{eqnarray*}
By (G2), the inner integral can be taken over 
$K - K = \{z - y: z \in K,\ y \in K\}$, and the sup-norm of the convolution 
is bounded by $\Vert Z(s,\cdot) \Vert^2_{L^2(\IR^d)}$, so this is
\[
  \leq E\left(\int_0^T ds\, \Vert Z(s,\cdot) \Vert^2_{L^2(\IR^d)} \Gamma (K - K)\right)
   = \int_0^T ds\, \Vert Z(s,\cdot) \Vert^2_{L^2(\IR^d)}\, \Gamma (K - K) < \infty,
\]
by (\ref{2.1}) and the fact that 
$s \mapsto E(\Vert Z(s,\cdot) \Vert^2_{L^2(\IR^d)})$ is continuous by (G2). 
Therefore, $v_{g,Z}(x)$ will be well-defined provided we show that 
$(s,y,\omega) \mapsto g(s,x-y) Z(s,y;\omega)$ is predictable, or 
equivalently, that $(s,y,\omega) \mapsto Z(s,y;\omega)$ is predictable.

   For this, set $t^n_j = j T 2^{-n}$ and
\[
   Z_n(s,y) = \sum_{j=0}^{2^n-1} Z(t^n_j, x)\, 1_{]t^n_j, t^n_{j+1}]}(s).\]
Observe that
\begin{eqnarray*}
   \Vert Z_n \Vert_+^2 &=& E\left(\sum_{j=0}^{2^n-1} \int_{t^n_j}^{t^n_{j+1}} ds 
         \int_{\IR^d} \Gamma(dx)\, (\vert Z(t^n_j,\cdot)\vert \ast 
         \vert \tilde Z(t^n_j,\cdot)\vert)(-x)\right) \\
        &\leq& T 2^{-n} \sum_{j=0}^{2^n-1} E(\Vert Z(t^n_j,\cdot)\Vert^2_{L^2(\IR^d)})\, \Gamma (K - K) \\
        &<& \infty.
\end{eqnarray*}
Therefore, $Z_n \in {\cal{P}}_+$, since this process, which is adapted, 
continuous in $x$ and left-continuous in $s$, is clearly predictable.  Further,
\begin{eqnarray*}
  && E\left(\int_0^T ds \int_{\IR^d} \Gamma(dx)\, (\vert Z(s,\cdot) - Z_n(s,\cdot)\vert) \ast 
       (\vert \tilde Z(s,\cdot) - \tilde Z_n(s, \cdot)\vert)(-x)\right) \\
 &&\qquad\qquad  \leq \int_0^T ds\, E(\Vert Z(s,\cdot) - Z_n(s,\cdot)\Vert^2_{L^2(\IR^d)}) 
    \, \Gamma (K - K).
\end{eqnarray*}
The integrand converges to $0$ and is uniformly bounded over $[0,T]$ by 
(G2), so this expression converges to $0$ as $n \to \infty$.  Therefore, $Z$ 
is predictable.

   Finally, we prove (\ref{2.3}).  Clearly, $E(\Vert v_{g, Z} \Vert^2_{L^2(\IR^d)})$ is equal to
\[
   E\left(\int_{\IR^d} dx\, \left(\int_0^T \int_{\IR^d} g(s, x-y) Z(s, y) M(ds, dy)\right)^2\right).
\]
Since the covariation measure of $M$ is $Q$, this equals
\begin{equation}\label{2.4}
   E\left(\int_{\IR^d} dx \int_0^T ds \int_{\IR^d} \Gamma(dz) \int_{\IR^d} dy 
   \, g(s, x-y) Z (s, y) g(s, x-z-y) Z(s, z+y)\right).
\end{equation}
The inner integral is equal to 
$(g(s, x-\cdot) Z(s, \cdot)) \ast (\tilde{g}(s, x - \cdot) \tilde{Z}(s, \cdot))(-z)$, 
and since this function belongs to $\IS(\IR^d)$ by (G1) and (G2), (\ref{2.4}) equals 
\[
   E\left(\int_{\IR^d} dx \int_0^T ds \int_{\IR^d} \mu (d \eta)\, 
   \vert \F(g(s, x - \cdot) Z(s, \cdot)) (\eta) \vert^2\right),
\]
by (\ref{2.2}).  Because the Fourier transform takes products to convolutions,
\[
   \F(g(s, x - \cdot) Z(s, \cdot))(\eta) = \int_{\IR^d} d \xi^\prime 
   e^{i \xi^\prime \cdot x} \F g(s, \cdot)(-\xi^\prime) \F Z (s, \cdot)(\eta - \xi^\prime),
\]
so, by Plancherel's theorem,
\begin{equation}\label{2.5}
   \int_{\IR^d} dx\, \vert \F(g(s, x-\cdot) Z(s, \cdot))(\eta) \vert^2 = \int_{\IR^d} 
    d \xi^\prime\, \vert \F g(s, \cdot)(-\xi^\prime) \F Z(s, \cdot) (\eta - \xi^\prime)\vert^2.
\end{equation}
The minus can be changed to plus, and using the change of variables 
$\xi = \eta + \xi^\prime$ ($\eta$ fixed), we find that (\ref{2.3}) holds. 
\hfill $\clubsuit$
\vskip 16pt

\begin{remark} An alternative expression for $I_{g,Z}$ is
\[
   I_{g,Z} = E\left(\int_0^T ds \int_{I\!\!R^d} \mu(d \eta)\, 
   \Vert g(s, \cdot) \ast (\chi_\eta(\cdot)Z(s, \cdot)) \Vert_{L^2(\IR^d)}^2\right),
\]
where $\chi_\eta(x) = e^{i \eta \cdot x}$.  Indeed, notice that (\ref{2.5}) is equal to
\[
   \int d \xi^\prime\, \vert \F g(s, \cdot)(\xi^\prime) \F(\chi_\eta(\cdot)Z(s, \cdot))(\xi^\prime)\vert^2,
\]
which, by Plancherel's theorem, is equal to
\[
   \Vert g(s, \cdot) \ast (\chi_\eta(\cdot)Z(s, \cdot)) \Vert^2_{L^2(\IR^d)}.
\]
\end{remark}

   Fix $(s, x) \mapsto g(s, x)$ such that (G1) holds.  Consider the further hypotheses:
\begin{itemize}
\item[{\bf (G4)}] $\int_0^T ds\, \sup_\xi \int_{\IR^d} \mu(d\eta)\, 
\vert \F g(s,\cdot)(\xi - \eta)\vert^2 < \infty$,
\item[{\bf (G5)}] For $0 \leq s \leq T$, $Z(s,\cdot) \in L^2(\IR^d)$ a.s., 
$Z(s,\cdot)$ is $\F_s$-measurable, and $s \mapsto Z(s,\cdot)$ is 
mean-square continuous from $[0,T]$ into $L^2(\IR^d)$.
\end{itemize}

   Fix $g$ such that (G1) and (G4) hold.  Set
\[
   \IP = \{Z \mbox{ : (G5) holds}\}.
\]
Define a norm $\Vert \cdot \Vert_g$ on $\IP$ by
\[
   \Vert Z \Vert_g^2 = I_{g,Z}.
\]
We observe that by (G4) and (G5) (and Plancherel's theorem), 
$I_{g,Z} \leq \tilde{I}_{g,Z} < \infty$, where
\[
   \tilde{I}_{g,Z} = \int_0^T ds\, E(\Vert Z(s,\cdot)\Vert^2_{L^2(\IR^d)})
      \left(\sup_\xi \int_{\IR^d} \mu(d\eta)\, 
      \vert \F g(s,\cdot)(\xi - \eta)\vert^2\right).
\]
Let
\[
   \E = \{Z \in \IP: \mbox{ (G2) holds}\}.
\]
By Lemma \ref{rlem1}, $Z \mapsto v_{g,Z}$ defines an isometry from 
$(\E,\Vert \cdot \Vert_g)$ into $L^2(\Omega \times \IR^d, dP \times dx)$. 
Therefore, this isometry extends to the closure of 
$(\E,\Vert \cdot \Vert_g)$ in $\IP$, which we now identify.

\begin{lemma} $\IP$ is contained in the closure of $(\E,\Vert \cdot \Vert_g)$.
\label{rlem3}
\end{lemma}

\proof Fix $\psi \in C_0^\infty(\IR^d)$ such that $\psi \geq 0$, the 
support of $\psi$ is contained in the unit ball of $\IR^d$ and 
$\int_{\IR^d} \psi(x)\, dx = 1$.  For $n \geq 1$, set 
\[
   \psi_n(x) = n^d \psi(nx).
\]
Then $\psi_n \to \delta_0$ in $\IS(\IR^d)$ and 
$\F \psi_n(\xi) = \F \psi(\xi/n)$, therefore 
$\vert {\cal{F}} \psi_n(\cdot)\vert$ is bounded by 1.

Fix $Z \in {\cal{P}}$, and show that $Z$ belongs to the completion of 
$\E$ in $\Vert \cdot \Vert_g$.  Set 
\[
   Z_n(s,x) = Z(s, x) 1_{[-n,n]^d}(x)\quad \mbox{ and }\quad Z_{n,m}(s, \cdot) = Z_n(s, \cdot) \ast \psi_m.
\]
We first show that $Z_{n,m} \in \E$, that is, (G2) holds for $Z_{n,m}$. 
Clearly, $Z_{n,m}(s, \cdot) \in C_0^\infty(\IR^d)$, $Z_{n,m}(s, \cdot)$ is 
$\F_s$-measurable by (G5), and there is a compact set 
$K_{n,m} \subset \IR^d$ such that supp~$Z_{n,m}(s, \cdot) \subset K$, for 
$0 \leq s \leq T$.  Further,
\[
  \Vert Z_{n,m}(t,\cdot) - Z_{n,m}(s,\cdot) \Vert^2_{L^2(\IR^d)} \leq \Vert 
  Z_n(t,\cdot) - Z_n(s,\cdot) \Vert^2_{L^2(\IR^d)} \leq \Vert Z(t,\cdot) - 
  Z(s,\cdot) \Vert^2_{L^2(\IR^d)},
\]
so $s \mapsto Z_{n,m}(s, \cdot)$ is mean-square continuous by (G5). 
Therefore, $Z_{n,m} \in \E$.

We now show that for $n$ fixed, $\Vert Z_n - Z_{n,m}\Vert_g \to 0$ as 
$m \to \infty$.  Clearly,
\[
   I_{g, Z_n-Z_{n, m}} = \int_0^T ds \int_{\IR^d} d \xi\, E(\vert \F Z_n(s, 
   \cdot)\vert^2) \ \vert 1 - \F \psi_m(\xi)\vert^2 \int_{\IR^d} \mu(d 
   \eta)\, \vert \F g(s, \cdot) (\xi - \eta)\vert^2.
\]
Because $\vert 1 - \F \psi_m(\xi)\vert^2 \leq 4$ and
\begin{eqnarray*}
  I_{g, Z_n} &\leq& \int_0^T ds\, E\left(\int_{\IR^d} d\xi\, \vert \F 
  	Z_n(s,\cdot)\vert^2\right)
      \left(\sup_\xi \int_{\IR^d} \mu(d \eta)\, \vert \F g(s, \xi -n)\vert^2 \right) \\
   &=& \tilde{I}_{g, Z_n} \leq \tilde{I}_{g, Z} < \infty,
\end{eqnarray*}
we can apply the Dominated Convergence Theorem to see that for $n$ fixed,
\[
   \lim_{m \to \infty} \Vert Z_n - Z_{n,m}\Vert_g 
   = \lim_{m \to \infty} \sqrt{I_{g, Z_n-Z_{n, m}}} = 0.
\]
Therefore, $Z_n$ belongs to the completion of $\E$ in $\Vert\cdot \Vert_g$. 
We now show that $\Vert Z-Z_n \Vert_g \to 0$ as $n \to \infty$.  Clearly,
\begin{eqnarray*}
   \Vert Z-Z_n \Vert_g^2 &=& I_{g, Z-Z_n} \leq \tilde{I}_{g, Z-Z_n} \\
    &=& \int_0^T ds\, 
    E\left(\Vert (Z - Z_n)(s,\cdot) \Vert^2_{L^2(\IR^d)}\right) 
    \left(\sup_\xi \int_{\IR^d} \mu(d \eta)\, \vert \F g(s, \xi -n)\vert^2 \right).
\end{eqnarray*}
Because 
\[
   \Vert Z-Z_n \Vert^2_{L^2(\IR^d)} \leq (\Vert Z \Vert_{L^2(\IR^d)} 
   + \Vert Z_n \Vert_{L^2(\IR^d)})^2 \leq 4 \Vert Z \Vert_{L^2(\IR^d)}^2,
\]
and $\tilde{I}_{g, Z} < \infty$, the Dominated Convergence Theorem implies that
\[
   \lim_{n \to \infty} \Vert Z - Z_n\Vert_g = 0,
\]
and therefore $Z$ belongs to the completion of $\E$ in $\Vert\cdot \Vert_g$. 
Lemma \ref{rlem3} is proved.
\hfill $\clubsuit$
\vskip 16pt

\begin{remark} Lemma \ref{rlem3} allows us to define the stochastic 
integral $v_{g,Z} = g \cdot M^Z$ provided $g$ satisfies (G1) and (G4), and 
$Z$ satisfies (G5).  The key property of this stochastic integral is that
\[
   E\left(\Vert v_{g,Z} \Vert_{L^2(\IR^d)}^2\right) = I_{g,Z}.
\]
\label{rrem4}
\end{remark}

We now proceed with a further extension of this stochastic integral, by 
extending the map $g \mapsto v_{g,Z}$ to a more general class of $g$.

Fix $Z \in \IP$.  Given a function $s \mapsto G(s) \in \IS'(\IR^d)$, 
consider the two properties:
\begin{itemize}
\item[{\bf (G6)}] For all $s \geq 0$, $\F G(s)$ is a function and 
\[
  \int_0^T ds\, \sup_\xi \int_{\IR^d} \mu(d \eta)\, \vert \F G(s, \xi - \eta)\vert^2 < \infty.
\]
\item[{\bf (G7)}] For all $\psi \in C_0^\infty(\IR^d)$, 
$\sup_{0 \leq s \leq T} G(s) \ast \psi$ is bounded on $\IR^d$.
\end{itemize}
Set 
\[
   \G = \{s \mapsto G(s): \mbox{ (G6) and (G7) hold}\},
\]
and
\[
   \IH = \{s \mapsto G(s): G(s) \in C^\infty(\IR^d) \mbox{ and (G1)and (G4) hold}\}.
\]
Clearly, $\IH \subset \G$.  For $G \in \G$, set
\[
   \Vert G \Vert_Z = \sqrt{I_{G,Z}}.
\]
Notice that $I_{G,Z} \leq \tilde{I}_{G,Z} < \infty$ by (G5) and (G6).  By 
Remark \ref{rrem4}, the map $G \mapsto v_{G,Z}$ is an isometry from 
$(\IH, \Vert \cdot \Vert_Z)$ into $L^2(\Omega \times \IR^d, dP \times dx)$. 
Therefore, this isometry extends to the closure of 
$(\IH, \Vert \cdot \Vert_Z)$ in $\G$.

\begin{lemma} $\G$ is contained in the closure of $(\IH, \Vert \cdot \Vert_Z)$.
\label{rlem5}
\end{lemma}

\proof Fix $s \mapsto G(s)$ in $\G$.  Let $\psi_n$ be as in the proof of 
Lemma \ref{rlem3}.  Set 
\[
  G_n(s, \cdot) = G(s) \ast \psi_n(\cdot).
\] 
Then $G_n(s, \cdot) \in C^\infty(\IR^d)$ by \cite{sch66}, 
Chap.VI, Thm.11 p.166.  By (G6), $\F G_n(s,\cdot) = \F G(s) \cdot \F\psi_n$ 
is a function, and so by (G7), (G1) holds for $G_n$.  Because 
$\vert \F \psi_n\vert \leq 1$, (G4) holds for $G_n$ because it holds for $G$ 
by (G6).  Therefore, $G_n \in \IH$.

Observe that 
\begin{eqnarray*}
   \Vert G - G_n\Vert_Z^2 &=& I_{G-G_n,Z} \\
    &=& \int^T_0 ds \int_{\IR^d} d \xi\, 
    E(\vert \F Z(s, \cdot)(\xi)\vert^2 \int_{\IR^d} \mu(d \eta)\, 
    \vert \F G(s,\cdot) (\xi- \eta) \vert^2 \vert 1 - \F \psi_n (\xi-\eta) \vert^2
\end{eqnarray*}
The last factor is bounded by $4$, has limit $0$ as $n \to \infty$, and 
$I_{G,Z} < \infty$, so the Dominated Convergence Theorem implies that
\[
   \lim_{n \to \infty} \Vert G - G_n\Vert_Z =0.
\]
This proves the lemma.
\hfill $\clubsuit$
\vskip 16pt

\begin{thm} Fix $Z$ such that (G5) holds, and $s \mapsto G(s)$ such that 
(G6) and (G7) hold.  Then the stochastic integral $v_{G,Z} = G \cdot M^Z$ is 
well-defined, with the isometry property
\[
   E\left(\Vert v_{G,Z} \Vert _{L^2(\IR^d)}^2\right) = I_{G,Z}.
\]
\label{rthm6}
\end{thm}

It is natural to use the notation
\[
   v_{G,Z} = \int^T_0 ds \int_{\IR^d} G(s,\cdot - y) Z(s,y) M(ds,dy),
\]
and we shall do this in the sequel.
\vskip 12pt

\noindent{\sc Proof of Theorem \ref{rthm6}.} The statement is an immediate 
consequence of Lemma \ref{rlem5}.
\hfill $\clubsuit$
\vskip 16pt

\begin{remark} Fix a deterministic function $\psi \in L^2(\IR^d)$ and set
\[
   X_t = \left\langle \psi,\ \int_0^t\int_{\IR^d} G(s,\cdot-y) Z(s,y)\,
           M(ds,dy)\right\rangle_{L^2(\IR^d)}.
\]
It is not difficult to check that $(X_t,\ 0 \leq t \leq T)$ is a 
(real-valued) martingale.
\end{remark}

\section{Examples}\label{sec3}
\setcounter{equation}{0}
In this section, we give a class of examples to which Theorem \ref{rthm6} 
applies.  Fix an integer $k \geq 1$ and let $G$ be the Green's function of 
the p.d.e.
\begin{equation}\label{3.1}
   \frac{\partial^2 u}{\partial t^2} + (-1)^k \Delta^{(k)} u = 0.
\end{equation}
As in \cite{dal99}, Section 3, $\F G(t)(\xi)$ is easily computed, and one finds
\[
   \F G(t)(\xi) = \frac{\sin(t\vert \xi \vert^k)}{\vert \xi \vert^k}.
\]
According to \cite{dal99}, Theorem 11 (see also Remark 12 in that paper), 
the linear s.p.d.e
\begin{equation}\label{3.2}
   \frac{\partial^2u}{\partial t^2} + (-1)^k \Delta^{(k)}u = {\dot F}(t, x)
\end{equation}
with vanishing initial conditions has a process solution if and only if
\[
\int_0^T ds \int_{I\!\!R^d} \mu (d \xi)\, \vert \F G (s) \vert^2 < \infty,
\]
or equivalently,
\begin{equation}\label{3.3}
   \int_{\IR^d} \mu (d \xi)\, \frac{1}{(1+\vert \xi \vert^2)^k} < \infty.
\end{equation}
It is therefore natural to assume this condition in order to study 
non-linear forms of (\ref{3.2}).

In order to be able to use Theorem \ref{rthm6}, we need the following fact.

\begin{lemma} Suppose (\ref{3.3}) holds.  Then the Green's function $G$ of 
equation (\ref{3.1}) satisfies conditions (G6) and (G7).
\label{rlem6}
\end{lemma}

\proof We begin with (G7).  For $\psi \in C_0^\infty (\IR^d)$,
\[
\begin{array}{ll}
   \Vert G(s) \ast \psi \Vert_{L^\infty(\IR^d)} &\leq \Vert \F(G(s) 
   \ast \psi)\Vert_{L^1(\IR^d)} = \displaystyle\int_{\IR^d} 
   \frac{\vert \sin(s \vert \xi \vert^k) \vert}{\vert \xi \vert^k} \vert \F\psi(\xi)\vert\, d \xi\\
&\leq s \displaystyle\int_{\IR^d} \vert \F \psi (\xi)\vert\, d \xi < \infty,
\end{array}
\]
so (G7) holds.

Turning to (G6), we first show that
\begin{equation}\label{3.4}
   \langle \chi_\xi G_{d,k} , \Gamma \rangle 
   = \left\langle \frac{1}{(1+ \vert \xi - \cdot \vert^2)^k}, \mu \right\rangle,
\end{equation}
where
\[
   0 \leq G_{d,k}(x) = \frac{1}{\gamma(k)} \int_0^\infty e^{-u} u^{k-1} p(u,x)\,
   du,
\]
$\gamma(\cdot)$ is Euler's Gamma function and $p(u, x)$ is the density of a 
$N(0, uI)-$ random vector (see \cite{san00}, Section 5).  In 
particular,
\[
   \F G_{d,k} (\xi) = \frac{1}{(1+\vert \xi \vert^2)^k},
\]
and it is shown in \cite{kar00} and \cite{san00} that 
\begin{equation}\label{3.5}
   \langle G_{d,k}, \Gamma \rangle = \langle (1+\vert \cdot \vert^2)^{-k}, \mu \rangle,
\end{equation}
and the right-hand side is finite by (\ref{3.3}).  However, the proofs in
\cite{san00} and \cite{kar00} use
monotone convergence, which is not applicable in presence of the 
oscillating function $\chi_\xi$.  As in \cite{kar00}, because 
$e^{-t \vert \cdot \vert^2}$ has rapid decrease,
\[
   \left\langle \frac{e^{-t \vert \cdot \vert^2}}{(1+ \vert \xi - \cdot \vert^2)^k}, 
   \mu \right\rangle = \left\langle \F \left(\frac{e^{-t \vert \cdot \vert^2}}
   {(1 + \vert \xi - \cdot \vert^2)^k}\right), \mu \right\rangle 
   = \langle p(t, \cdot) \ast (\chi_\xi G_{d,k}), \Gamma \rangle.
\]
Notice that $G_{d,k} \geq 0$, and so
\[
   \vert p(t, \cdot) \ast (\chi_\xi G_{d,k})\vert \leq p(t, \cdot) \ast G_{d,k} 
   \leq e^T G_{d, k}
\]
by formula (5.5) in \cite{san00}, so we can use monotone convergence in the first 
equality below and the Dominated Convergence Theorem in the third equality 
below to conclude that
\[
\begin{array}{lll}
   \left\langle \frac{1}{(1+ \vert \xi - \cdot \vert^2)^k}, \mu \right\rangle 
   &=& \lim_{t \downarrow 0} \left\langle \frac{e^{-t \vert \cdot \vert^2}}
   {(1+ \vert \xi - \cdot \vert^2)^k} , \mu\right\rangle 
   = \lim_{t \downarrow 0} \langle p(t, \cdot) \ast (\chi_\xi G_{d,k})), \Gamma \rangle \\
   && \\
   &=& \langle \lim_{t \downarrow 0} (p(t, \cdot) \ast (\chi_\xi G_{d,k})), \Gamma \rangle
   = \langle \chi_\xi G_{d, k}, \Gamma \rangle,
\end{array}
\]
which proves (\ref{3.4}).  Because $G_{d,k} \geq 0$,
\begin{equation}\label{3.6}
   \sup_\xi\, \langle \chi_\xi G_{d,k}, \Gamma\rangle \leq \langle G_{d,k}, \Gamma\rangle < \infty 
\end{equation}
by (\ref{3.5}) and (\ref{3.3}).  The lemma is proved.
\hfill $\clubsuit$
\vskip 16pt

\section{A non-linear s.p.d.e}\label{sec4}
\setcounter{equation}{0}
Let $\alpha: \IR \to \IR$ be a Lipschitz function such that $\alpha(0) = 0$, 
so that there is a constant $K > 0$ such that for $u, u_1, u_2 \in \IR$,
\begin{equation}\label{4.1}
   \vert \alpha(u)\vert \leq K \vert u \vert \qquad 
   \mbox{and} \qquad \vert \alpha(u_1) - \alpha(u_2) \vert \leq K \vert u_1 - u_2 \vert. 
\end{equation}
Examples of such functions are $\alpha(u) = u$, $\alpha(u) = \sin(u)$, or $\alpha(u) = 1-e^{-u}$. 

   Consider the non-linear s.p.d.e.
\begin{equation}\label{4.2}
   \frac{\partial^2 }{\partial t^2}u(t,x) + (-1)^k \Delta^{(k)}u(t,x) = \alpha(u(t, x)) F(t, x),
\end{equation}
\[
   u(0, x) = v_0(x),\qquad \frac{\partial}{\partial t} u(0, x) = \tilde{v}_0(x) 
\]
where $v_0 \in L^2(\IR^d)$ and $\tilde{v}_0 \in H^{-k}(\IR^d)$, the Sobolev  
space of distributions such that
\[
   \Vert \tilde{v}_0 \Vert^2_{H^{-k}(\IR^d)} 
   \stackrel{\mbox{\scriptsize def}}{=} \int_{\IR^d} d\xi\, 
   \frac{1}{(1+\vert \xi \vert^2)^k} \vert \F \tilde{v}_0(\xi)\vert^2 < \infty.
\]
We say that a process $(u(t, \cdot),\ 0 \leq t \leq T)$ with values in 
$L^2(\IR^d)$ is a solution of (\ref{4.2}) if, for all $t \geq 0$, a.s.,
\begin{equation}\label{4.3}
   u(t, \cdot) = \frac{d}{dt} G(t) \ast v_0 + G(t) \ast \tilde{v}_0 
   + \int_0^t \int_{\IR^d} G(t-s, \cdot -y) \alpha(u(s, y)) M(ds, dy),
\end{equation}
where $G$ is the Green's function of (\ref{3.1}).  The third term is 
interpreted as the stochastic integral from Theorem \ref{rthm6}, so 
$(u(s, \cdot))$ must be adapted and mean-square continuous from $[0, T]$ 
into $L^2(\IR^d)$.

\begin{thm} Suppose that (\ref{3.3}) holds.  Then equation (\ref{4.2}) 
has a unique solution $(u(t, \cdot),\ 0 \leq t \leq T)$.  This solution is 
adapted and mean-square continuous.
\label{rthm7}
\end{thm}

\proof We will follow a standard Picard iteration scheme.  Set
\[
   u_0(t, \cdot) = \frac{d}{dt} G(t) \ast v_0 + G(t) \ast \tilde{v}_0.
\]
Notice that $v_0(t, \cdot) \in L^2(\IR^d)$.  Indeed,
\begin{eqnarray}\label{4.4}
   \left\Vert \frac{d}{dt} G(t) \ast v_0\right\Vert_{L^2(\IR^d)} 
   &=& \left\Vert \F \frac{d}{dt} G(t) \cdot \F v_0\right\Vert_{L^2(\IR^d)} 
   = \int_{\IR^d} \sin^2(t\vert \xi \vert^k) \vert \F v_0 (\xi)\vert^2 d\xi\\
   &\leq& \Vert v_0\Vert_{L^2(\IR^d)},\nonumber
\end{eqnarray}
and one checks similarly that 
$\Vert G(t) \ast \tilde{v}_0 \Vert_{L^2(\IR^d)} \leq \Vert \tilde{v} \Vert_{H^{-k}}$. 
Further, $t \mapsto u_0(t, \cdot)$ from $[0, T]$ into $L^2(\IR^d)$ is continuous. Indeed, 
\[
   \lim_{t \to s} \left\Vert \frac{d}{dt}G(t) 
   \ast v_0 - \frac{d}{dt} G(s) \ast v_0 \right\Vert_{L^2(\IR^d} = 0, 
\]
as is easily seen by proceeding as in (\ref{4.4}) and using dominated 
convergence. Similarly, 
\[
   \lim_{t \to s} \Vert  G(t) \ast \tilde{v}_0 - G(s) \ast \tilde{v}_0\Vert = 0.
\]
 
   For $n \geq 0$, assume now by induction that we have defined an adapted and 
mean-square continuous process $(u_n(s, \cdot),\ 0 \leq s \leq T)$ with 
values in $L^2(\IR^d),$ and define
\begin{equation}\label{4.5}
   u_{n+1}(t, \cdot) = u_0 (t, \cdot) + v_{n+1}(t, \cdot),
\end{equation}
where
\begin{equation}\label{4.6}
   v_{n+1} (t, \cdot) = \int_0^t \int_{\IR^d} G(t-s, \cdot -y) \alpha(v_n(s, g)) M(ds, dy). 
\end{equation}
We note that $(\alpha(u_n(s, \cdot)),\ 0 \leq s \leq T)$ is adapted and 
mean-square continuous, because by (\ref{4.1}),
\[
   \Vert \alpha(u_n(s, \cdot)) - \alpha(u_n(t, \cdot))\Vert_{L^2(\IR^d)} 
   \leq K \Vert u_n(s, \cdot) - u_n(t, \cdot) \Vert_{L^2(\IR^d)},
\]
so the stochastic integral in (\ref{4.6}) is well-defined by Lemma 
\ref{rlem6} and Theorem \ref{rthm6}.

Set
\begin{equation}\label{4.7}
   J(s) = \sup_\xi \int_{\IR^d} \mu(d \eta)\, \vert \F G(s, \cdot)(\xi - \eta)\vert^2.
\end{equation}
By (\ref{3.3}), (\ref{3.5}) and (\ref{3.6}), 
$\sup_{0 \leq s \leq T} J(s)$ is bounded by some $C < \infty$, so by 
Theorem \ref{rthm6} and using (\ref{4.1}),
\begin{eqnarray}\nonumber
   E(\Vert u_{n+1}(t, \cdot) \Vert^2_{L^2(\IR^d)}) &\leq& 2 \Vert u_0 (t, \cdot) \Vert^2_{L^2(\IR^d)} 
        + 2 \int_0^t ds\, E(\Vert \alpha(u_n(s, \cdot))\Vert^2_{L^2(\IR^d)}) J(t-s)\\
   &\leq& 2 \Vert u_0 (t, \cdot)\Vert_{L^2(\IR^d)}^2 + 2 K C \int_0^t ds\, E(\Vert u_n(s, \cdot)\Vert^2_{L^2(\IR^d)}).
\label{4.8}
\end{eqnarray}
Therefore, $u_{n+1}(t, \cdot)$ takes its values in $L^2(\IR^d)$.  By Lemma 
\ref{rlem7} below, $(u_{n+1}(t, \cdot),\ 0 \leq t \leq T)$ is mean-square 
continuous and this process is adapted, so the sequence $(u_n,\ n \in \IN)$ 
is well-defined.  By Gronwall's lemma, we have in fact
\[
   \sup_{0 \leq t \leq T} \sup_{n \in \IN} E(\Vert u_n(t, \cdot) \Vert^2_{L^2(\IR^d)}) < \infty.
\]
We now show that the sequence $(u_n(t, \cdot),\ n \geq 0)$ converges.  Let
\[
   M_n(t) = E(\Vert u_{n+1} (t, \cdot) - u_n(t, \cdot) \Vert^2_{L^2(\IR^d)}).
\]
Using the Lipschitz property of $\alpha(\cdot), (\ref{4.5})$ and 
$(\ref{4.6})$, we see that
\[
   M_n(t) \leq K C \int_0^t ds\, M_{n-1}(s).
\]
Because $\sup_{0 \leq s \leq T} M_0(s) < \infty$, Gronwall's lemma implies that
\[
   \sum_{n=0}^\infty M_n(t)^{1/2} < \infty.
\]
In particular, $(u_n(t, \cdot),\ n \in \IN)$ converges in 
$L^2(\Omega \times \IR^d, dP \times dx)$, uniformly in $t \in [0, T]$, to a 
limit $u(t, \cdot)$.  Because each $u_n$ is mean-square continuous and the 
convergence is uniform in $t$, $(u(t, \cdot),\ 0 \leq t \leq T)$ is also 
mean-square continuous, and is clearly adapted.  This process is easily seen 
to satisfy (\ref{4.3}), and uniqueness is checked by a standard argument.  
\hfill $\clubsuit$
\vskip 16pt

The following lemma was used in the proof of Theorem \ref{rthm7}.

\begin{lemma} Each of the processes $(u_n(t, \cdot),\ 0 \leq t \leq T)$ 
defined in (\ref{4.5}) is mean-square continuous.
\label{rlem7}
\end{lemma}

\proof Fix $n \geq 0$.  It was shown in the proof of Theorem \ref{rthm7} 
that $t \mapsto u_0(t, \cdot)$ is mean-square continuous, so we establish 
this property for $t \mapsto v_{n+1}(t, \cdot)$, defined in $(\ref{4.6})$. 
Observe that for $h > 0$,
\[
   E(\Vert v_{n+1}(t+h, \cdot) - v_{n+1} (t, \cdot) \Vert^2_{L^2(\IR^d)}) \leq 2 (I_1 + I_2),
\]
where
\begin{eqnarray*}
   I_1 &=& E\left(\left\Vert \int_t^{t+h} \int_{\IR^d} G(t+h-s, \cdot - y) 
   \alpha(u_n(s, y)) M(ds, dy)\right\Vert^2_{L^2(\IR^d)}\right),\\
   I_2 &=& E\left(\left\Vert \int_0^t \int_{\IR^d} 
   (G(t+h-s, \cdot -y) - G(t-s, \cdot -y)) \alpha(u_n(s, y)) 
              M(ds dy) \right\Vert^2_{L^2(\IR^d)}\right).
\end{eqnarray*}
Clearly,
\[
   I_1 \leq K^2 \int_t^{t+h} ds\, E( \Vert u_n(s, \cdot) \Vert^2_{L^2(\IR^d)}) J(t+h-s),
\]
while
\[
   I_2 = \int_0^t ds \int_{\IR^d} d \xi\, \vert \F(\alpha(u_n(s, \cdot)))(\xi)\vert^2 
   \int_{\IR^d} \mu(d \eta) \left(\frac{\sin(t+h-s) \vert \xi - \eta \vert^k) 
   - \sin((t-s) \vert \xi - \eta \vert^k)}{\vert \xi - \eta \vert^k}\right)^2.
\]
The squared ratio is no greater than
\[
   4 \left( \frac{\sin(h \vert \xi - \eta \vert^k)}
   {\vert \xi - \eta \vert^k}\right)^2 
   \leq \frac{C}{(1 + \vert \xi - \eta \vert^2)^k}.
\]
It follows that $I_2$ converges to $0$ as $h \to 0$, by the dominated 
convergence theorem, and $I_1$ converges to $0$ because the integrand is 
bounded.  This proves that $t \mapsto v_{n+1}(t, \cdot)$ is mean-square 
right-continuous, and left-continuity is proved in the same way. 
\hfill $\clubsuit$
\vskip 16pt

\section{The wave equation in weighted $L^2$-spaces}\label{sec5}
\setcounter{equation}{0}
In the case of the wave equation (set $k=1$ in (\ref{4.2})), we can 
consider a more general class of non-linearities $\alpha(\cdot)$ than in 
the previous section.  This is because of the compact support property of 
the Green's function of the wave equation. 

More generally, in this section, we fix $T > 0$ and consider a function 
$s \mapsto G(s) \in {\cal{S}}'(\IR^d)$ that satisfies (G6), (G7) and, in 
addition,
\begin{itemize}
\item[{\bf (G8)}] There is $R > 0$ such that for $0 \leq s \leq T$, 
supp~$G(s) \subset B(0, R)$.
\end{itemize}

Fix $K > d$ and let $\theta : \IR^d \to \IR$ be a smooth function for which 
there are constants $0 < c < C$ such that
\[
   c(1 \wedge \vert x \vert^{-K} ) \leq \theta(x) \leq C(1 \wedge \vert x \vert^{-K}).
\]
The weighted $L^2$-space $L^2_\theta$ is the set of measurable 
$f: \IR^d \to \IR$ such that $\Vert f \Vert_{L_\theta^2} < \infty$, where
\[
   \Vert f \Vert_{L_\theta^2}^2 = \int_{\IR^d} f^2(x) \theta(x)\, dx.
\]
Let $H_n = \{x \in \IR^d: n R \leq \vert x \vert < (n+1)R\}$, set 
\[
   \Vert f \Vert_{L^2(H_n)} = \left(\int_{H_n} f^2(x)\, dx\right)^{1/2},
\]
and observe that there are positive constants, which we again denote $c$ and $C$, 
such that
\begin{equation}\label{5.1}
   c \sum_{n=0}^\infty n^{-K} \Vert f \Vert^2_{L^2(H_n)} 
   \leq \Vert f \Vert^2_{L_\theta^2} 
   \leq C \sum_{n=0}^\infty n^{-K} \Vert f \Vert^2_{L^2(H_n)}.
\end{equation}
For a process $(Z(s, \cdot),\ 0 \leq s \leq T)$, consider the following hypothesis:
\begin{itemize}
\item[{\bf (G9)}] For $0 \leq s \leq T$, $Z(s, \cdot) \in L_\theta^2$ a.s., 
$Z(s, \cdot)$ is $\F_s$-measurable, and $s \mapsto Z(s, \cdot)$ is mean-square continuous from $[0, T]$ into $L^2_\theta$.
\end{itemize}
Set
\[
\begin{array}{lll}
   \E_\theta = \{Z: &\mbox{ (G9) holds, and there is } K \subset \IR^d \mbox{ compact such}\\
& \mbox{that for } 0 \leq s \leq T, \mbox{ supp} \, Z(s, \cdot) \subset K \}.
\end{array}
\]
Notice that for $Z \in \E_\theta$, $Z(s, \cdot) \in L^2(\IR^d)$ because 
$\theta(\cdot)$ is bounded below on $K$ by a positive constant, and for the 
same reason, $s \mapsto Z(s, \cdot)$ is mean-square continuous from $[0, T]$ 
into $L^2(\IR^d)$.  Therefore,
\[
   v_{G,Z} = \int_0^T \int_{\IR^d} G(s, \cdot - y) Z(s, y) M(ds, dy)
\]
is well-defined by Theorem \ref{rthm6}.

\begin{lemma} For $G$ as above and $Z \in \E_\theta$, $v_{G, Z} \in L^2_\theta$ a.s. and
\[
E(\Vert v_{G, Z} \Vert^2_{L^2_\theta}) \leq \int_0^T ds\, \Vert Z(s, \cdot) \Vert^2_{L_\theta^2} J(s),
\]
where $J(s)$ is defined in $(\ref{4.7})$.
\label{rlem8}
\end{lemma}

\proof We assume for simplicity that $R = 1$ and $K$ is the unit ball in 
$\IR^d$.  Set $D_0 = H_0 \cup H_1$ and, for $n \geq 1$, set 
$D_n = H_{n-1} \cup H_n \cup H_{n+1}$ and 
$Z_n(s, \cdot) = Z(s, \cdot) 1_{D_n}(\cdot)$.  By (\ref{5.1}), then (G8),
\begin{eqnarray*}
  \Vert v_{G,Z} \Vert^2_{L_\theta^2} &\leq& \sum_{n=0}^\infty n^{-K} \Vert v_{G,Z} \Vert^2_{L^2(H_n)} 
             = \sum_{n=0}^\infty n^{-K} \Vert v_{G,Z_n} \Vert^2_{L^2(H_n)}\\
   &\leq& \sum_{n=0}^\infty n^{-K} \Vert v_{G,Z_n} \Vert^2_{L^2(\IR^d)}.
\end{eqnarray*}
Therefore, by Theorem \ref{rthm6},
\begin{eqnarray*}
   E(\Vert v_{G,Z} \Vert^2_{L^2_\theta} 
   &\leq& \sum_{n=0}^\infty n^{-K} \int_0^T ds\, E(\Vert Z_n(s, \cdot) \Vert^2_{L^2(I\!\!R^d)}) J(s)\\
   &=& \sum_{n=0}^\infty n^{-K} \int_0^T ds\, E(\Vert Z(s, \cdot) \Vert^2_{L^2(D_n)}) J(s)\\
&\leq& C \int_0^T ds\, \sum_{n=0}^\infty n^{-K}\, E(\Vert Z(s, \cdot) \Vert^2_{L^2(H_n)}) J(s)\\
&\leq& C \int_0^T ds\, E( \Vert Z(s, \cdot) \Vert^2_{L^2_\theta})J(s).
\end{eqnarray*}
This proves the lemma. 
\hfill $\clubsuit$
\vskip 16pt

For a process $(Z(s, \cdot))$ satisfying (G9), let 
$\Vert Z \Vert_\theta =(I^\theta_{G, Z})^\half$, where
\[
   I_{G, Z}^\theta = \int_0^T ds\, \Vert Z(s, \cdot) \Vert_{L_\theta^2}^2 J(s).
\]
Because $s \mapsto \Vert Z(s, \cdot) \Vert_{L_\theta^2}^2$ is bounded, 
$I^\theta_{G,Z} < \infty$ provided (G6) holds.
Therefore, $\Vert Z \Vert_\theta$ defines a norm, and by Lemma \ref{rlem8}, 
$Z \mapsto v_{G,Z}$ from $\E_\theta$ into 
$L^2(\Omega \times \IR^d, dP \times \theta(x)dx)$ is continuous.  Therefore 
this map extends to the closure of $\E_\theta$ for 
$\Vert \cdot \Vert_\theta$, which we now identify.

\begin{thm} Consider a function $s \mapsto G(s) \in {\cal{S}}'(\IR^d)$ such 
that (G6), (G7) and (G8) hold.  Let $(Z(s, \cdot),\ 0 \leq s \leq T)$ be an 
adapted process with values in $L_\theta^2$ that is mean-square continuous 
from $[0, T]$ into $L^2_\theta$.  Then $Z$ is in the closure of $\E_\theta$ 
for $\Vert \cdot \Vert_\theta$, and so the stochastic integral $v_{G, Z}$ 
is well-defined, and
\begin{equation}\label{5.2}
   E(\Vert v_{G,Z} \Vert^2_{L^2_\theta}) \leq I^\theta_{G, Z}.
\end{equation}
\label{rthm8}
\end{thm}

\proof Set $Z_n(s, \cdot) = Z(s, \cdot) 1_{[-n,n]^d}(\cdot)$.  Then $(Z_n)$ 
satisfies (G9) and belongs to $\E_\theta$.  Because 
$\Vert Z_n(s, \cdot)\Vert_{L^2_\theta}\leq \Vert Z(s, \cdot) \Vert_{L^2_\theta}$ 
and $I_{G,Z}^\theta< \infty$, the dominated convergence theorem implies 
that $\lim_{n \to \infty} \Vert Z-Z_n\Vert_\theta = 0$, so $Z$ is in the 
closure of $\E_\theta$ for $\Vert \cdot \Vert_\theta$, and $(\ref{5.2})$ 
holds by Lemma \ref{rlem8}.
\hfill $\clubsuit$
\vskip 16pt

   We now use this result to obtain a solution to the following stochastic wave equation:
\begin{equation}\label{5.3}
   \frac{\partial^2}{\partial t^2}u(t,x) - \Delta u(t,x) = \alpha(u(t, x)) {\dot F}(t, y),
\end{equation}
\[
   u(0, x) = v_0(x),\qquad \frac{\partial u}{\partial t} (0, x) = \tilde{v}_0(x).
\]
We say that a process $(u(t, \cdot),\ 0 \leq t \leq T)$ with values in 
$L^2_\theta$ is a solution of $(\ref{5.3})$ if $(u(t, \cdot))$ is 
adapted, $t \mapsto u(t, \cdot)$ is mean-square continuous from $[0, T]$ 
into $L^2_\theta$ and
\begin{equation}\label{5.4}
   u(t, \cdot) = \frac{d}{dt} G(t) \ast v_0 + G(t) \ast \tilde{v}_0 
   + \int_0^t \int_{\IR^d} G(t-s, \cdot -y) \alpha(u(s, y)) M(ds, dy),
\end{equation}
where $G$ is the Green's function of the wave equation.  In particular, 
$\F G(s) (\xi) = \vert \xi \vert^{-1} \sin(t \vert \xi \vert)$ and (G6), 
(G7) and (G8) hold provided $(\ref{3.3})$ holds with $k=1$.  Therefore, the 
stochastic integral in $(\ref{5.4})$ is well-defined by Theorem \ref{rthm8}.

\begin{thm} Suppose
\[
   \int_{\IR^d} \mu(d\xi)\, \frac{1}{1+\vert \xi\vert^2} < \infty,
\]
$v_0 \in L^2(\IR^d)$, $\tilde{v}_0 \in H^{-1}(\IR^d)$, and $\alpha(\cdot)$ 
is a globally Lipschitz function.  Then (\ref{5.3}) has a unique solution 
in $L^2_\theta$.
\label{rthm9}
\end{thm}

\proof The proof follows that of Theorem \ref{rthm7}, so we only point out 
the changes relative to the proof of that theorem.  Because $\alpha(\cdot)$ 
is globally Lipschitz, there is $K > 0$ such that for $u, u_1, u_2 \in \IR$,
\[
   \vert \alpha (u) \vert \leq K(1+\vert u\vert)\qquad 
   \mbox{ and }\qquad \vert \alpha(u_1) - \alpha(u_2) \vert \leq K \vert u_1-u_2\vert.
\]
Using the first of these inequalities, $(\ref{4.8})$ is replaced by
\[
   E(\Vert u_{n+1} (t, \cdot) \Vert^2_{L^2_\theta}) \leq 2 \Vert u_0(t, \cdot) 
   \Vert_{L^2(\IR^d)} + 2 KC \int_0^t ds\, (1 + E(\Vert u_n(s, \cdot) \Vert^2_{L_\theta^2})).
\]
Therefore $u_{n+1}(t, \cdot)$ takes its values in $L^2_\theta$.  The 
remainder of the proof is unchanged, except that 
$\Vert \cdot \Vert_{L^2(\IR^d)}$ must be replaced by 
$\Vert \cdot \Vert_{L_\theta^2}$.  This proves Theorem \ref{rthm9}.
\hfill $\clubsuit$


\vskip 16pt

\begin{thebibliography}{99}
\bibliographystyle{alpha}
\bibitem{ada75} Adams, R.A. {\em Sobolev Spaces.} Pure and Applied Mathematics,
Vol. 65, Academic Press, New York-London, 1975.

\bibitem{cn88a} Carmona, R. and Nualart, D. Random nonlinear wave equations: propagation of singularities
{\em Annals Probab.} 16 (1988), 730-751.

\bibitem{cn88b} Carmona, R. and Nualart, D. Random nonlinear wave equations: smoothness of the solutions
{\em Probab. Theory Related Fields} 79 (1988), 469--508.

\bibitem{dal99} Dalang, R.C. Extending the martingale measure stochastic integral witha
pplications to spatially homogeneous s.p.d.e's. {\em Electron. J. Probab.} 4
(1999), 29pp.

\bibitem{dz92} Da Prato, G. and Zabczyk, J. {\em Stochastic Equations in
Infinite Dimensions}. Encyclopedia of mathematics and its applications 44. 
Cambridge University Press, Cambridge, New York, 1992.

\bibitem{df98} Dalang, R.C. and Frangos, N. E. The stochastic wave equation in
two spatial dimensions. {\em Annals Probab.} 26 (1998) 187-212.

\bibitem{kar00} Karkzewska, A. and Zabczyk, J. Stochastic PDEs with
function-valued solutions. Preprint 33, Scuola Normale Superiore di Pisa
(1997).

\bibitem{kr81} Krylov, N.V. and Rozovskii, B.L. Stochastic evolution systems. 
{\em J. Soviet Math.} 16 (1981), 1233-1276.

\bibitem{kr82} Krylov, N.V. and Rozovskii, B.L. Stochastic partial differential equations
and diffusion processes. {\em Russian Math. Surveys} 37 (1982), 81-105.

\bibitem{Lev01} L\'ev\`eque, O. {\em Hyperbolic stochastic partial differential equations driven by
boundary noises.} Ph.D.~thesis, no.2452, Ecole Polytechnique F\'ed\'erale de 
Lausanne, Switzerland (2001).

\bibitem{mill99} Millet, A. and Sanz-Sol\'e, M. A stochastic wave equation in two space dimensions: smoothness of
the law. {\em Annals Probab.} 27 (1999), 803-844.

\bibitem{mue97} Mueller, C. Long time existence for the wave equation with a noise
term. {\em Annals Probab.} 25 (1997), 133-152.

\bibitem{or98} Oberguggenberger, M. and Russo, F. Nonlinear stochastic wave
equations. {\em Integral Transform. Spec. Funct.} 6 (1998), 71-83.

\bibitem{par72} Pardoux, E. Sur des \'equations aux d\'eriv\'ees partielles stochastiques
monotones. {\em C. R. Acad. Sci. Paris} S\'er. A-B 275 (1972), A101-A103.

\bibitem{par75} Pardoux, E. Equations aux d\'eriv\'ees partielles stochastiques de
type monotone. S\'eminaire sur les \'Equations aux D\'eriv\'ees Partielles
(1974--1975), III, Exp. No. 2 (1975), p.10.

\bibitem{par77} Pardoux, E. Characterization of the density of the conditional law in the
filtering of a diffusion with boundary. In: Recent developments in statistics 
(Proc. European Meeting Statisticians, Grenoble, 1976). North Holland,
Amsterdam (1977), 559-565.

\bibitem{pes02} Peszat, S. The Cauchy problem for a nonlinear stochastic wave equation 
in any dimension. To appear in {\em J. Evol. Equ.} (2002).

\bibitem{pz00} Peszat, S. and Zabczyk, J. Nonlinear stochastic wave and heat
equations. {\em Probab. Theory Related Fields} 116 (2000), 421-443.

\bibitem{san00} Sanz-Sol\'e, M. and Sarr\`a, M. Path properties of a class of martingale measures with applications
to spde's. In: {\em Stochastic processes, physics and geometry: new interplays, I
(Leipzig, 1999).} (Gestesy, F., Holden, H., Jost, J., Paycha, S., R\"ockner, 
M.~and Scarlatti, S., eds). CMS Conf.~Proc., Amer.~Math.~Soc., Providence, RI
(2000), 345-355.

\bibitem{sch66} Schwartz, L. {\em Th\'eorie des distributions.} Hermann, Paris,
1966.

\bibitem{ste70} Stein, E.M. {\em Singular Integrals and Differentiability Properties
of Functions.}  Princeton University Press, Princeton, New Jersey, 1970.

\bibitem{wal86} Walsh, J.B. An Introduction to Stochastic Partial Differential
Equations. In: {\em Ecole d'Et\'e de Probabilit\'es de Saint-Flour, 
XIV-1984}, Lecture Notes in Mathematics 1180. Springer-Verlag, Berlin, Heidelberg, 
New York (1986), 265-439.

\end{thebibliography}
\end{document}